\documentclass[11pt]{article}
\usepackage{latexsym,amssymb}
\usepackage{a4wide}
\usepackage[T1]{fontenc}
\usepackage{fourier}
\usepackage{theorem}
\usepackage[svgnames]{pstricks}
\usepackage{fancyhdr}
\usepackage{hyperref}
\usepackage{sectsty}

\hypersetup{colorlinks=true,linkcolor=violet,citecolor=purple}
\allsectionsfont{\scshape}
\sectionfont{\large\scshape}
\subsectionfont{\normalsize\scshape}
\psset{arrowscale=1.5,dotscale=1.5}

\DeclareMathAlphabet{\mathcal}{OMS}{cmsy}{m}{n}

\def\headcolour{\color{DarkSeaGreen}}
\newtheorem{theorem}{Theorem}

\newtheorem{proposition}[theorem]{Proposition}

{\theorembodyfont{\rmfamily}
 
 \newtheorem{example}[theorem]{Example}
 
 \newtheorem{remark}[theorem]{Remark}
 
}
\newenvironment{proof}{\addvspace\baselineskip\noindent{\sc
    Proof}\hspace{0.5ex}}{\hspace*{\fill} $\Box$\par\addvspace\baselineskip}

\def\defn#1{{\bfseries\itshape #1}}
\def\restr#1{\lower4pt\hbox{$\rule[-2pt]{0.5pt}{11pt}\,\scriptstyle#1$}}

\def\RR{\mathbb{R}}
\def\NN{\mathbb{N}}
\def\ZZ{\mathbb{Z}}
\def\one{\mathbb{1}}
\def\mm{{m}}

\def\CS{\mathcal{C}} 
\def\ISG{\mathcal{I}} 
\def\cT{\mathcal{T}} 

\def\Fix{\mathop\mathsf{Fix}\nolimits}
\def\tr{\mathop\mathsf{tr}\nolimits}

\def\eps{\varepsilon}

\def\SO{\mathsf{SO}}
\def\OO{\mathsf{O}}
\def\Tet{\mathbb{T}}
\def\Oct{\mathbb{O}}
\def\Icos{\mathbb{I}}
\def\C{\mathsf{C}} 
\def\D{\mathsf{D}} 

\begin{document}
\thispagestyle{empty}

\noindent{\huge\textbf{Existence of symmetric central configurations}}

\bigskip
\bigskip

{\color{gray}
\noindent{\Large\bf James Montaldi}

\medskip

\noindent{\bf University of Manchester}

\medskip

\noindent March 4, 2015

\bigskip
}

{\small
\noindent\hrulefill 

\smallskip

\noindent{\large\scshape\color{gray} Abstract}

\medskip

\noindent 
Central configurations have been of great interest over many years, with the earliest examples due to Euler and Lagrange.  There are numerous results in the literature  demonstrating the existence of central configurations with specific symmetry properties, using slightly different techniques in each.  The aim here is to describe a uniform approach by adapting to the symmetric case the well-known variational argument showing the existence of central configurations.  The principal conclusion is that there is a central configuration for every possible symmetry type, and for any symmetric choice of masses.  Finally the same argument is applied to the class of balanced configurations introduced by Albouy and Chenciner.
\medskip


\noindent\hrulefill 
}

\bigskip

\section*{Introduction}

In the $n$-body problem, central configurations allow particularly simple motions.  If the particles are released from a central configuration with zero initial velocity, the configuration will collapse to the centre of mass while maintaining the same shape up to rescaling. If they are given other particular initial velocities, each particle will follow an elliptical Kepler orbit, and the shape formed by the configuration will remain constant up to rescaling and rotation. They also occur as limiting configurations of parabolic motions \cite{M-V09} and (partial) collisions \cite{F-T04}. R.~Moeckel has written a recent survey on the subject \cite{Moeckel14}.

We consider the set of central configurations in $\RR^d$.  Of course, the most interesting cases are $d=2$ and $d=3$, but nothing is lost by considering general dimensions.  Over the past few decades, many papers have been written demonstrating the existence of central configurations with various different symmetries, for example  \cite{CeLl89, CoLl10, CoLl13, L-S06, Yu-Zh12, Zh-Ch13} (and references therein).   The aim of this paper is to describe a uniform proof of all these existence results, using well-known arguments for the existence of symmetric solutions to variational problems.  
The main result is the following.
\medskip

\noindent\textbf{Theorem}\quad \emph{Given any symmetric configuration of $n$ bodies in $\RR^d$ and a corresponding symmetric distribution of masses, there is at least one central configuration of that symmetry type and with the given masses.} 

\medskip

\begin{figure} 
\psset{dotsize=3pt,linewidth=0.5pt,unit=1.2}
\centering
\begin{pspicture}(-2,-2)(3,2)  
{  \psset{linecolor=lightgray}
   \psline(-2,0)(2,0)
   \psline(-1,-1.732)(1,1.732)
   \psline(1,-1.732)(-1,1.732)
}
{  \psset{unit=1.6} 
   \psdots[linecolor=Blue](-1,0)(0.5,0.866)(0.5,-0.866)
   \psline[linestyle=dashed,linecolor=gray](-1,0)(0.5,0.866)(0.5,-0.866)(-1,0)
}
{\psset{unit=0.7}
   \psdots[linecolor=Crimson](-1,0)(0.5,0.866)(0.5,-0.866)
   \psline[linestyle=dashed,linecolor=gray](-1,0)(0.5,0.866)(0.5,-0.866)(-1,0)
}
{\psset{unit=0.4} 
\psdots[linecolor=LimeGreen](3., 1.)(-.6340, 3.098)(-2.366, 2.098)(-2.366, -2.098)(-.6340, -3.098)(3, -1)(3,1)
\psline[linestyle=dashed,linecolor=gray](3., 1.)(-.6340, 3.098)(-2.366, 2.098)(-2.366, -2.098)(-.6340, -3.098)(3, -1)(3,1)
}
\end{pspicture}
\begin{minipage}{0.75\textwidth}
\caption{A configuration with triangular ($\D_3$) symmetry, consisting of 12 points forming 3 orbits: two equilateral triangles and one semiregular hexagon. The theorem guarantees the existence of a configuration of this form, where the relative sizes will depend on the relative masses of each of the 3 orbits.}
\label{fig:D_3 configuration}
\end{minipage}
\end{figure}
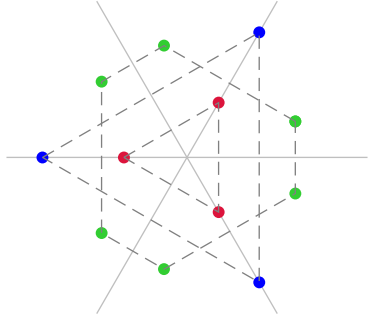

We state a more precise result as Theorem\,\ref{thm:main} below, after defining what is meant by symmetry type (or Burnside type), and a refinement using connected components.  An example of a symmetric configuration with triangular ($D_3$) symmetry is illustrated in Figure\,\ref{fig:D_3 configuration}: {in order to be a central configuration}, the relative sizes of the three orbits will depend on the relative masses of each.  The proof of the theorem uses the well-known variational approach to existence of central configurations, adapted to the symmetric setting, and details are given in Section\,\ref{sec:proof} below.  

In Section~\ref{sec:examples} we give a few examples in 2 and 3 dimensions. We show for example the existence of nested and staggered (or dual) platonic solids, as well as (nested) cubeoctahedron and icosidodecahedron configurations, and discuss (in Example\,\ref{eg:Archimedean}) why other Archimedean configurations are not likely to be central. In Section~\ref{sec:topology} we briefly describe the topological aspect of this problem.  The final short section illustrates how the same techniques can be applied to \emph{balanced configurations}, an extension of the idea of central configuration due to Albouy and Chenciner \cite{A-C98}, and we make a few observations about the relation between symmetric central and balanced configurations.

\section{Symmetric configurations}
\label{sec:symmetry}

A configuration of $n$ particles in $\RR^d$ is simply a set of $n$ points in $\RR^d$, each having a mass.  It is usual to order these points, so that the configuration is given as a point in $(\RR^d)^n=\RR^{nd}$. However, in order to avoid introducing permutation groups when we consider symmetric configurations, we wish to avoid ordering the points. We therefore describe a configuration of $n$ particles as a set $C:=\{x_1,\dots,x_n\}\subset\RR^d$ together with a function $\mm:C\to \RR^+$, assigning to each point in $C$ the mass of the particle at that point.  We denote the set of such configurations $(C,m)$ of $n$ particles in $\RR^d$ by $\CS\left(\RR^d,n\right)$, or simply by $\CS$. 

The group $\OO(d)$ of orthogonal transformations consists of rotations and rotation-reflections (or, improper rotations) in $\RR^d$, the former having determinant equal to 1, the latter to $-1$ (recall that an orthogonal transformation of determinant $-1$ is the product of a reflection and a rotation).  This group acts on the space of configurations in the natural way: let $(C,m)$ be a configuration with $C=\left\{x_1,\dots,x_n\right\}$ and let $g\in\OO(d)$, then $g\cdot(C,m) = (g\cdot C,\,g\cdot m)$ where
\begin{equation} \label{eq:action on C}
g\cdot C = \left\{gx_1,\dots,gx_n\right\},
\end{equation}
and for the mass function, $g\cdot\mm:g\cdot C\to \RR^+$ is defined by $(g\cdot\mm)(g\cdot x) =\mm(x)$.  That is, $g\cdot\mm=\mm\circ g^{-1}$.

Now consider a finite subgroup $G$ of the orthogonal group $\OO(d)$.  A configuration $(C,m)$ is a \defn{symmetric configuration} if the group $G$ leaves the set invariant: $g\cdot C=C$ (it will usually permute the points within the set), and moreover it preserves the masses, so that $g\cdot\mm=\mm$. In particular, this requires that the points $x$ and $gx$ in the configuration have the same mass; we call this an invariant mass distribution.  Since all our arguments and results are independent of the mass provided it is an invariant mass distribution, we may in places ignore the mass function, and concentrate just on the configuration of points.

For a given finite subgroup $G$ of $\OO(d)$, the symmetric configurations therefore form the subset $\CS^G:=\Fix(G,\CS)$ of $\CS$.   Let $C$ be a symmetric configuration (with invariant mass function $m$). If $x\in C$ then so is $gx$, and therefore so is the \defn{orbit} of $x$, which is the set of images of $x$ under the elements of $G$:
$$G\cdot x \ = \ \left\{gx\in\RR^d \mid g\in G\right\}.$$

If $G$ acts on a finite set, then the set can be partitioned into a disjoint union of orbits (as in Figure~\ref{fig:D_3 configuration} where the 12 points form 3 orbits).  However, different orbits may have different `geometry', and this is made precise by the  \emph{orbit type} of an orbit defined as follows.  The \defn{isotropy subgroup} $G_x$ of a point $x$ is the subgroup of $G$ consisting of those transformations fixing $x$:
$$G_x \ = \ \left\{g\in G\mid gx = x\right\}.$$
In particular, if $x=0$ then $G_x=G$. It is a simple exercise to show that if $y=gx$ then $G_y=gG_xg^{-1}$; that is the isotropy subgroups of two points in the same orbit are conjugate.  Thus to each orbit is associated a \emph{conjugacy class} of subgroups of $G$, called the \defn{orbit type} of the orbit. For a subgroup $H$ of $G$, one denotes the conjugacy class containing $H$ by $(H)$, and for  $x\in\RR^d$, the orbit type of $x$ is therefore $(G_x)$. The number of points in an orbit of type $(H)$ is equal to $|G|/|H|$ (where $|H|$ is the order of a group $H$).

Notice in particular that if an orbit has type $(H)$ say, then at least one of the points $x$ of the orbit has isotropy subgroup $G_x=H$ and thus lies in the fixed point subspace
$$\Fix(H,\,\RR^d) = \left\{x\in\RR^d \mid hx=x,\;\;\forall h\in H\right\},$$
which is a linear subspace of $\RR^d$. 

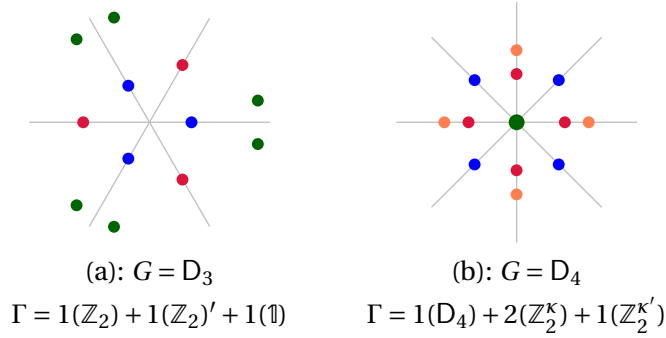
\begin{figure} 
\psset{dotsize=3pt,linewidth=0.5pt,unit=0.8}
\centering
\begin{pspicture}(-2,-3.5)(3,2)  
 {\psset{linecolor=lightgray}
   \psline(-2,0)(2,0)
   \psline(-1,-1.732)(1,1.732)
   \psline(1,-1.732)(-1,1.732)
 }
{\psset{unit=-0.7}    \psdots[linecolor=Blue](-1,0)(0.5,0.866)(0.5,-0.866)}
{\psset{unit=1.1} \psdots[linecolor=Crimson](-1,0)(0.5,0.866)(0.5,-0.866)}
{\psset{unit=1.2} \psdots[linecolor=DarkGreen](1.5,0.3)(1.5,-0.3)(-1.01,1.15)(-1.01,-1.15)(-0.490,1.45)(-0.490,-1.45) }
 \rput(0,-2.5){(a): $G=\D_3$}
 \rput(0,-3.2){$\Gamma=1(\ZZ_2)+1(\ZZ_2)'+1(\one)$}
\end{pspicture}
\begin{pspicture}(-3,-3.5)(2,2)  
 {\psset{linecolor=lightgray}
   \psline(-2,0)(2,0)
   \psline(0,-2)(0,2) 
   \psline(1.414,1.414)(-1.414,-1.414)
   \psline(1.414,-1.414)(-1.414,1.414)
 }
 \psdot[dotsize=4pt,linecolor=DarkGreen](0,0)
{\psset{unit=0.8} \psdots[linecolor=Crimson](1,0)(-1,0)(0,1)(0,-1)}
{\psset{unit=1.2} \psdots[linecolor=Coral](1,0)(-1,0)(0,1)(0,-1)}  
{\psset{unit=0.7} \psdots[linecolor=Blue](1,1)(1,-1)(-1,1)(-1,-1)} 
 \rput(0,-2.5){(b): $G=\D_4$}
 \rput(0,-3.2){$\Gamma=1(\D_4)+2(\ZZ^{\kappa}_2)+1(\ZZ^{\kappa'}_2)$}
\end{pspicture}
\begin{minipage}{0.75\textwidth}
\caption{$\D_3$ and $\D_4$ symmetric configurations in the plane, taking connected components into account for $\D_3$.}
\label{fig:D_n configurations}
\end{minipage}
\end{figure}

As a simple example consider the dihedral subgroup $D_3$ of $\OO(2)$; this is the symmetry group of the equilateral triangle in the plane. See Figures\,\ref{fig:D_3 configuration} and \ref{fig:D_n configurations}(a). There are, at this point in the discussion, three types of orbit: the origin with orbit type $(D_3)$, an orbit of type $(\ZZ_2)$ consisting of 3 points forming an equilateral triangle  (each vertex of the triangle is fixed by a reflection) and finally a `generic' orbit of type $(\one)$ consisting of 6 points forming a semiregular hexagon all with trivial isotropy.  In this way we can write a general $D_3$-symmetric configuration as an integer combination of orbit types (so many orbits of each orbit type): we write
$$\Gamma = \eps(\D_3) + a(\ZZ_2) + b(\one),$$
where $\eps\in\{0,1\}$ (since the only point with isotropy $\D_3$ in the plane is the origin, so there can be at most one such orbit), while $a,b\in\NN=\{0,1,2,\dots\}$.  A similar discussion applies to $\D_4$, see Figure\,\ref{fig:D_n configurations}(b), but note that there are two non-conjugate reflections in $\D_4$, here denoted $\kappa$ and $\kappa'$.

Extending this example to $\RR^3$, we let $\D_3$ act as before on the $(x,y)$-coordinates, and $\ZZ_2$ act by reflection in the $(x,y)$-plane, so $\tau(x,y,z)=(x,y,-z)$ (the Schoenflies notation for this subgroup of $\OO(3)$ is $\D_{3h}$). There are now a total of 6 orbit types:  the three considered above in the plane $z=0$, but now with isotropy type enhanced by $\ZZ_2^\tau$, so for example the orbit type of the 3-point orbit is ($\ZZ_2^\kappa\times\ZZ_2^\tau)$, where $\kappa$ is a reflection in $\D_3$, and three new ones which are, firstly a pair of opposite points on the $z$-axis at $(0,0,\pm z)$ for some $z\neq0$, with isotropy type $(\D_3)$, secondly orbits of 6 points forming a triangular prism, which has isotropy type $(\ZZ_2^\kappa)$, and finally the `generic' orbit consisting of 12 points arranged at the vertices of a semiregular hexagonal prism, with isotropy type $(\one)$. Thus a general symmetric configuration with symmetry $\D_{3h}$ is of the form
$$\Gamma = \eps(\D_3\times\ZZ_2^\tau) + a(\ZZ_2^\tau\times\ZZ_2^\kappa) + b(\ZZ_2^\tau) + c(\D_3) + d(\ZZ_2^\kappa) +e(\one),$$
with again $\eps\in\{0,1\}$ and $a,b,c,d,e\in\NN$.

This idea of writing a $G$-invariant set as an integer combination of orbit types goes back to Burnside \cite{Burnside} in the early days of group theory, so we call this the \defn{Burnside type} of a symmetric configuration.  Many details about properties of Burnside types set can be found in Kerber \cite{Kerber}.  

In order to treat equilateral triangles and their `duals' as distinct types, we need to refine the Burnside type to what we call the \emph{topological} Burnside type. This is most easily illustrated with the simple $\D_3$ example shown in Figure\,\ref{fig:D_n configurations}(a).  While the two triangles (red and blue) have the same orbit type, they cannot be continuously deformed one into the other whilst maintaining that orbit type, and so belong to different connected components of the set of orbits with orbit type $(\ZZ_2)$.  We denote these as $(\ZZ_2)$ and  $(\ZZ_2)'$. A similar phenomenon occurs in $\RR^3$ with tetrahedra and their duals.  In contrast, Figure\,\ref{fig:D_n configurations}(b) shows that the square and its dual have distinct isotropy types, here denoted $(\ZZ_2^\kappa)$ and  $(\ZZ_2^{\kappa'})$, where $\kappa$ is the reflection in the $x$-axis and $\kappa'$ the reflection in the diagonal $y=x$, and $\ZZ_2^\kappa$ denotes the group of order 2 generated by $\kappa$; for $\D_4$ the topology does not refine the Burnside type.

We therefore define the \defn{topological Burnside type} by distinguishing connected components of the set of orbits of type $(H)$ into connected components, writing them as $(H), (H)'$ etc, or more generally  $(H)^\alpha$ for $\alpha$ in some index set.

We are now in a position to state a more precise version of the theorem above. 

\begin{theorem} \label{thm:main}
Given any finite subgroup $G$ of\/ $\OO(d)$ and any topological Burnside type $\Gamma$ for $G$, there is at least one central configuration in each connected component of the set\/ $\CS(\Gamma)$.
\end{theorem}

The set $\CS(\Gamma)$ fails to be connected only if one of the fixed point spaces is 1-dimensional and the number of orbits of the corresponding type is greater than 1, for then reordering those points {may correspond} to different connected components.  See Remark\,\ref{rmk:ordered} and Section\,\ref{sec:topology} for more details.

The precise central configuration whose existence is given by the theorem will depend on the values of the masses of the particles (recall that for a symmetric configuration the mass distribution is invariant: that is, points in the same orbit have equal mass). From these existence theorems, under non-degeneracy conditions which for most mass distributions will be generic, one can apply the implicit function theorem to obtain central configurations with non-symmetric mass distributions, at least for nearby values of the masses, though the configurations will no longer be symmetric in general. Moreover when the central configurations are degenerate one expects to see bifurcations, some of which may break the symmetry, as for example in \cite{Moeckel-Simo}.  (There are also bifurcations which do not break the symmetry such as in \cite{L-S06}).

{
\begin{remark} \label{rmk:ordered}
We have described configurations of particles as \emph{unordered} sets of points, together with the mass of each one.  It is more traditional to describe configurations as \emph{ordered} collections of points $(x_1,x_2,\dots,x_n)$, with respective masses $(m_1,m_2,\dots,m_n)$.  There are two reasons for adopting our approach. Firstly, the results are independent of ordering, and it is artificial to introduce the (arbitrary) ordering. Secondly, the symmetry group $G$ no longer leaves the ordered configuration invariant, but it permutes the elements, so for each $g\in G$ there is a permutation $\sigma=\sigma(g)\in S_n$ such that $g\cdot x_i = x_{\sigma(i)}$.  The proof would require the action of $(g,\sigma)$ rather than just $g$.  The two approaches are in fact equivalent, and the Burnside type can still be defined using the `ordered' approach.  

One consequence of this approach is perhaps surprising. Consider for example configurations of two distinct particles in a line. If we ignore the mass, then a configuration is of the form $\{x,y\}$ with $x\neq y\in\RR$, and the set of such pairs is connected since, for example, $\{1,2\}=\{2,1\}$. On the other hand,  if we consider the set of ordered pairs $\{(x,y)\mid x\neq y\}$ then there are two connected components, one with $x<y$ and the other with $y<x$.  Now include the masses: the graph of $m$ is $\{(x,m(x)), \;(y,m(y))\}$.  If $m(x)=m(y)$ then there is again only one component, while if $m(x)\neq m(y)$ there are two: one with the larger mass on the right, the other with the larger mass on the left. Thus the topology of the configuration space as we define it depends on the masses, and this would not be the case if we consider ordered configurations. This might seem undesirable at first sight, but I claim it is completely natural: consider for example the collinear Euler relative equilibria for the 3-body problem. If the masses are distinct there are 3 such configurations (up to rescaling), depending which particle lies between the other two (recall that opposite orderings are equivalent under rotation in the plane). On the other hand, if they have the same mass, the 3 solutions are really the same, and there is just the one solution. 
\end{remark}
}

\section{Proof of the theorem} \label{sec:proof}
This type of theorem is usually presented by representing the configurations as ordered $n$-tuples and then using the permutation group acting by permuting the points as in \cite{LMR01,MS13}; this approach is needed particularly if collisions are involved, such as in \cite{St96}.  However, this is not necessary for our problem and here we proceed directly on the configurations as sets, as described above, which removes the need for introducing permutations. 

Without loss of generality, we restrict attention to configurations whose centre of mass is at the origin: $\sum_{x\in C}m(x)x=0$. Central configurations are determined by two functions defined for any configuration (with mass).  First the potential,
$$U(C) = \sum_{\{x,y\}\subset C}\frac{m(x)m(y)}{\|x-y\|},$$
where $m(x)$ is the mass of the particle at the point $x$ and the sum is over all unordered pairs of distinct points in the configuration $C$.  The other function is the total {moment of} inertia about the origin,
$$I(C)=\sum_{x\in C}m(x)\|x\|^2.$$
We can take as our definition the following, a configuration $C$ is a \defn{central configuration} if $C$ is a critical point of $U$ when restricted to a level set of $I$.  Since both functions are homogeneous (of degrees -1 and 2 respectively)  it follows that if a configuration  $C$ is central then so is the homothetic configuration $\lambda C = \{\lambda x\mid x\in C\}$ for any $\lambda\neq 0$, and consequently we can restrict attention to the level set $I=1$ for convenience.  See for example \cite{Moeckel14} for details. Let $\CS_1=\CS\cap\{I=1\}$.

Since both $U$ and $I$ depend only on the distances between the particles (and their masses), they are both invariant under the orthogonal group $\OO(d)$.  To prove the theorem, we use the so-called \emph{principle of symmetric criticality}, first established by Palais. But first we recall some basic facts about fixed point spaces.  If a group $G$ acts smoothly on a manifold $M$, then the set of points with symmetry $G$ is the fixed point set $M^G$, that is the subset
$$M^G:=\Fix(G,M)=\{x\in M\mid G\cdot x = x\}.$$
If the group is finite (or indeed compact) then if non-empty, $M^G$ is a union of closed submanifolds of $M$, possibly having components of different dimensions.

\medskip

\noindent \textbf{Principle of symmetric criticality} (Palais \cite{Palais}) 
\textit{Suppose a finite group $G$ acts smoothly on a manifold $M$ and suppose $f:M\to\RR$ is a smooth invariant function, and let $x\in M^G$.  Then $x$ is a critical point of $f$ if and only if it is a critical point of the restriction $f\restr{M^G}$.}

\medskip 

This is because at any symmetric point $x\in M^G$ the gradient $\nabla f(x)$ is tangent to $M^G$.

It follows from this principle that $C$ is a central configuration with symmetry $G$ and moment of inertia $I=1$ if and only if it is a critical point of the restriction of $U$ to the closed submanifold $\CS_1^G$ of $\CS_1$. 

{The manifold structure and topology of $\CS$ are defined simply by identifying it locally with the same configurations considered as ordered collections of points, regardless of which ordering is chosen, while the mass function is taken to be locally constant.  With this topology, one can show that each $\C(\Gamma)$ is both oen and closed in $\C^G$, and hence is a union of connected components of $\C^G$.}

Now $I:\CS\to\RR$ is a smooth non-singular invariant function so that $\CS_1=I^{-1}(1)$ is a smooth $G$-invariant submanifold, and hence the fixed point supspace $\CS_1^G$ is a smooth submanifold. 
Thus we have a smooth function $U:\CS_1^G\to \RR$, and we want to show it must have a critical point, for then the result follows by the principle of symmetric criticality.

Because of the form of $I$ (positive definite quadratic form),  if $C_j$ is a sequence of configurations in $\CS_1$ or in $\CS_1^G$ that doesn't contain a limit point, then the minimal distance between pairs of points must tend to zero.  Consequently, $U\to\infty$ on such a sequence. It follows that $U$ must attain a minimum somewhere on $\CS_1$, or $\CS_1^G$ respectively, and this minimum is the desired critical point. \hfill$\Box$

\begin{remark}
As is well-known, if $C$ is a central configuration in $\RR^d$ and $e>d$ then $C$ is also a central configuration in $\RR^{e}$, when embedded in $\RR^d\times\{0\}\subset \RR^{e}$.  To see this using a symmetry argument, write $\RR^{e} =\RR^d\times\RR^{e-d}$, and consider the 2-element subgroup of $\OO(e)$ generated by $g=\pmatrix{I&0\cr 0&-I}$. Then $\Fix(g,\RR^{e})=\RR^d$, and the result follows from the principle of symmetric criticality.
\end{remark}

\begin{remark}
We have been considering an ambient space of arbitrary dimension $d$.  If $d>3$ the relevance of the inverse square law is debatable, and for physical reasons should arguably be replaced by an inverse power $(d-1)$ law. The potential would then be of the form
$$U(C) = \sum_{\{x,y\}\subset C}\frac{m(x)m(y)}{\;\|x-y\|^{d-2}\;},$$
However, Theorem~\,\ref{thm:main}  only relies on the symmetry of the function $U$ and the fact that as the configuration approaches a collision, so $U\to\infty$.  It follows that the approach would also apply with this gravitational law, and indeed with any other potential depending only on the shape of the configuration and the masses, in a symmetric fashion, provided it tends to infinity near collisions. 
\end{remark}

\section{Examples}
\label{sec:examples}

{There are many statements in the literature of the existence of symmetric central configurations, and all can be deduced from the method described in this paper. In this section we describe a few of these.}

\subsection{Dimension 2}
Here it is straightforward to list the different types of symmetric configuration. The only finite subgroups of $\OO(2)$ are the cyclic groups $\C_k$ (of order $k$) and the dihedral groups $\D_k$ (of order $2k$).

If $G=\C_k$ ($k>1$), then there are two relevant Burnside types (and no topological refinement), namely $(\C_k)$ and $(\one)$. The corresponding fixed point subspaces are the origin and $\RR^2$. A symmetric configuration is then a set of $n=ak$ points, forming $a$ regular $k$-gons centred on the origin, together with possibly a point at the origin.  The Burnside type is 
\begin{equation}\label{eq:Ck Burnside}
\Gamma=\eps(\C_k)+a(\one)
\end{equation}

The group $G=\D_1$ consists of a reflection in (say) the $x$-axis. There are two (topological) Burnside types $(\D_1)$ (points on the axis) and ($\one$) (pairs of points, each a reflection of the other in the axis).  Therefore, given any Burnside type $\Gamma=a(\D_1)+b(\one)$, there is at least one central configuration of this type. 

Now suppose $G=\D_k$ ($k>1$). In this case there are four topological Burnside types, as described above for $k=3$ or 4.
\begin{equation} \label{eq:Dk Burnside}
\left\{\begin{array}{ll}
(\D_k),\; (\ZZ_2^\kappa),\;(\ZZ_2^\kappa)',\;(\one) & \mbox{for $k$ odd}\\[4pt]
(\D_k),\; (\ZZ_2^\kappa),\;(\ZZ_2^{\kappa'}),\;(\one) & \mbox{for $k$ even}
\end{array}\right.
\end{equation}
To treat the two cases together, denote these topological orbits as $(\D_k),(A),(B)$ and $(\one)$. Then a general symmetric configuration would have topological Burnside type
$$\Gamma=\eps(\D_k)+a(A)+b(B)+c(\one).$$
Geometrically, this would consist of $\eps$ points at the origin, $a$ regular nested $k$-gons, $b$ staggered (or twisted by $\pi/k$) regular $k$-gons and $c$ semiregular $2k$-gons, all centred at the origin.  A semiregular $2k$-gon is the orbit of a point in the complement of the axes of reflection {(this may in fact be a regular $2k$-gon, but that would not be a consequence of the $\D_k$-symmetry)}.  The study in \cite{L-S06} considers configurations with (topological) Burnside type $(\D_k)+(A)+(B)$. 

In each of these cases, the theorem guarantees the existence of central configurations with such symmetry, provided the masses have corresponding symmetry. We therefore recover and extend a result of Zhao and Chen \cite{Zh-Ch15} (their result corresponds to symmetry $\D_k$ and topological Burnside type $p(A)+g(B)$).

Since $\C_k<\D_k$ it follows that dihedral configurations also have cyclic symmetry.  However, to the best of my knowledge it is unknown whether there exist central configurations of equal mass with cyclic symmetry that do not in fact have dihedral symmetry.  On the other hand, by altering the masses of a configuration with dihedral symmetry, it is possible to produce one with precisely cyclic symmetry and no more.  For example, consider the Burnside type $(\one)$ with $G=\D_k$ (analogous to the semiregular hexagon in Figure~\ref{fig:D_3 configuration}). Now perturb the mass $m$ of alternate particles to a nearby value $m'$.  The perturbed central configuration will then not have any reflectional symmetry, but the rotation remains; it will therefore have symmetry $\C_k$ and be of Burnside type $2(\one)$.  I am grateful to Alain Albouy for this observation.

\subsection{Dimension 3}
This case is more complex, resulting in many more types of central configuration. A description of all the possible symmetry types is given in \cite{LMR01}, although some adaptation is needed as in that reference the action is restricted to the sphere: in particular the origin did not appear and nested polyhedra are not possible. 

\begin{example} \label{eg:C2h}  
Consider $G=\ZZ_2\times\ZZ_2$ with one generator $\tau$ acting by reflection in the $(x,y)$-plane, and the other $\rho$ by rotation by $\pi$ about the $z$-axis.  The Schoenflies notation is $\C_{2h}$. There are 4 Burnside types: $(C_{2h}),\, (\ZZ_2^\rho),\, (\ZZ_2^\tau)$ and $(\one)$. The Burnside type $\Gamma= \eps(C_{2h}) + a(\ZZ_2^\rho) + b(\ZZ_2^\tau)+c(\one)$ consists of $\eps$ points at the origin, $2a$ points symmetrically placed along the $z$-axis, $2b$ points in the $(x,y)$-plane placed symmetrically with respect to the origin and 4c points in space, placed in $G$-orbits (these each forming the vertices of a rectangle).  For any symmetric distribution of masses among these points, the theorem tells us that there is at least one central configuration with the points in such a configuration. 
\end{example}

\begin{example}\label{eg:Dnh}
Consider the subgroup $\D_{nh}<\OO(3)$ in the Schoenflies notation. As a group this is isomorphic to $\ZZ_2\times\D_n$, and is generated by the reflection $\tau$ in the $(x,y)$-plane (giving the $\ZZ_2$ factor) and the usual dihedral group acting on the $(x,y)$-plane and leaving the `vertical' $z$-axis fixed. Among the orbit types are the origin with orbit type  $(\D_{nh})$, the horizontal lines of reflection forming two components (as in 2 dimensions) with orbit type $(\ZZ_2^\tau\times\ZZ_2^\kappa)$ and $(\ZZ_2^\tau\times\ZZ_2^{\kappa'})$ or $(\ZZ_2^\tau\times\ZZ_2^\kappa)'$ (accordingly as $n$ is even or odd) giving orbits of regular $n$-gons and their duals, and the prisms with $n$-fold symmetry with the $2n$ vertices lying in the vertical planes of reflection and with $z\neq 0$ and orbit type $(\ZZ_2^\kappa)$ and $(\ZZ_2^{\kappa'})$ or $(\ZZ_2^{\kappa})'$ as above. Consider in particular the configurations with $3n$ points with Burnside type
$$\Gamma = \cases{1(\ZZ_2^\kappa)+1(\ZZ_2^\tau\times\ZZ_2^{\kappa'})& if $n$ is even\cr 
  1(\ZZ_2^\kappa)+1(\ZZ_2^\tau\times\ZZ_2^{\kappa})'& if $n$ is odd.}$$
This is chosen so that the $n$-gon in the plane $z=0$ is staggered (dual) relative to the polygons in the other horizontal planes.  The theorem then implies there must be a central configuration of this symmetry type, so proving the existence part of a conjecture of Corbera and Llibre \cite{CoLl13} on `double antiprisms'.  A similar result is available if the three $n$-gons are aligned rather than staggered. 
\end{example}

\begin{example}\label{eg:Platonic}
Consider the symmetry group $\Tet_d$ of the regular tetrahedron, which has order 24. There are 5 orbit types: $(\Tet_d)$ (the origin), $(S_3)$ (radial lines through the vertices of the tetrahedron or its dual), $(\ZZ_2\times\ZZ_2)$ (mid-points of the 6 edges, forming an octahedron), $(\ZZ_2)$ (other points on the edges, forming an orbit of 12 points) and $(\one)$ (generic points, orbits of 24 points).   The theorem tells us that for any non-negative integers $\eps,a,b,c,d,e$, there is a symmetric configuration of Burnside type
$$C = \eps(\Tet_d) + a(S_3)+b(S_3)'+c(\ZZ_2\times\ZZ_2)+d(\ZZ_2) + e(\one).$$
Here as usual $\eps\in\{0,1\}$ determines whether or not there is a point at the origin, $a$ is the number of nested tetrahedra and $b$ the number of nested dual tetrahedra, etc. 

Similar results apply to the other groups $\Oct_h$ and $\Icos_h$, from which we deduce a stronger form of the existence theorem of \cite{CoLl10} on nested Platonic solids---{in that paper they only show there exist masses for which such central configurations exist}.  

One can also deduce the existence of two types of Archimedean solid: the cubeoctahedron and the icosidodecahedron. 
The vertices of the cubeoctahedron lie at the mid-points of the cube (or of the octahedron) and has octahedral symmetry $\Oct_h$; it is uniquely determined by the orbit type  $(\ZZ_2\times\ZZ_2)$, a subgroup generated by reflections in two orthogonal planes, and this shows that this is also a central configuration.  The icosidodecahedron consists of 30 vertices placed at the mid-points of the edges of the dodecahedron (or of the icosahedron) and is similarly determined uniquely by the analogous orbit type $(\ZZ_2\times \ZZ_2)$, but now as a subgroup of the icosahedral group $\Icos_h$, and this shows it too is a central configuration.  Similarly, nested cubeoctahedra and nested icosidodecahedra also form central configurations. 
\end{example}

\begin{example}\label{eg:Archimedean}
On the other hand, other Archimedean solids do not (in all likelihood) form central configurations. This is because their symmetry group does not determine their shape.  For example, consider the family of truncated tetrahedra. These are obtained by shaving off the 4 vertices of a regular tetrahedron, replacing them with 3 vertices each and 4 new equilateral triangles as faces.  The original faces of the tetrahedron then become semiregular hexagons. As more is shaved off, the ratio between the lengths of the semiregular hexagons varies (increases say), and when the two lengths are equal, the hexagon is regular, and this truncated tetrahedron is an Archimedean solid.  Let $\rho>0$ denote the ratio of the sides of the semiregular hexagon. As $\rho\to0$ so the orbit tends to a tetrahedron, and as $\rho\to\infty$ the 12 vertices merge in pairs to form an octahedron. The Archimedean truncated tetrahedron of course corresponds to $\rho=1$. The theorem implies that there is at least one value of $\rho>0$ which forms a central configuration, and numerical calculations (using Maple) suggest this to be unique with value $\rho=0.855$ which does not correspond to the Archimedean shape (the edge between two semiregular hexagons being shorter than the edges of the equilateral triangles).

It is to be expected that a similar phenomenon happens for the other Archimedean shapes: namely that they fail to be central configurations, except of course the cubeoctahedron and icosidodecahedron discussed above. The distinguishing feature of these two particular shapes among the Archimedean ones is that they are \emph{edge regular}, which means that all the edges are equivalent under the symmetry group, while in the others there are two distinct types of edge, and the symmetry alone does not force them to be of equal length.
\end{example}

\begin{example} \label{eg:D2}
Consider the subgroup of $\SO(3)$ with 4 elements consisting of the identity, and the rotations by $\pi$ about each of the $x,y$ and $z$-axes. The Schoenflies notation is $\D_2$.  There are 5 types of orbit for this group.  Firstly the point at the origin which is fixed by the whole group. Secondly, a pair of opposite points on the $z$-axis, and these have isotropy equal to the subgroup of order 2 generated by the corresponding rotation $R_z$, third and fourth are the corresponding pairs of points on the $x$- and $y$-axes, and finally a generic orbit consisting of the 4 points $\{(x,y,z), (x,-y,-z),$  $(-x,y,-z),$ $(-x,-y,z)\}$ which are distinct provided at most one of the coordinates is 0, and which has trivial isotropy. Any symmetric configuration has Burnside type
$$\Gamma = \eps(G) + a(R_z) + b(R_y) + c(R_x) + e(\one).$$
($d$ is used for dimension!) 
In particular, with $\eps=0,a=k,b=p,c=\ell,e=0$ we reclaim the result of Jiang and Zhao \cite{J-Zh} using their notation. And of course we can let $\eps=1$ or $e>0$ to obtain a more general result.
\end{example}

\section{Topology}
\label{sec:topology}

Let $G$ be a given finite subgroup of $\OO(d)$ and let $\Gamma$ be a Burnside type for $G$. It is natural (and useful) to have a measure of the complexity of the corresponding set of configurations $\CS(\Gamma)$ using topological invariants.  This information can be used to find a lower bound on the number of central configurations of the given type, using Morse theory if all critical points are non-degenerate or more generally using Lusternik-Schnirelman category.  In this short section we give some indications of what this topology is. 

Consider the quotient space $X=\RR^d/G$, which is in general a singular space.  {The image of the points with orbit type $(H)$ is a subset of $X$, which we denote $X_H$, and which has a manifold structure such that together the $X_H$ form a stratification of $X$.} Note that if $H$ and $H'$ are conjugate then $X_H=X_{H'}$.  See for example \cite{DK00} for details on group actions and stratifications of their orbit space. 

Thus, if a particular symmetric configuration consists of $a$ orbits of type $(H)$, then it is determined by $a$ points in $X_H$.  Denote by $\ISG=\ISG(G)$ the collection of all conjugacy classes of isotropy subgroups of $G$, and by $\cT$ its refinement into topological classes (the $(H)^\alpha$ introduced earlier). For a conjugacy class $(H)\in \ISG$ we work with a representative $H$ and we write $a_H$ for the number of orbits of type $(H)$.

{Recall that a configuration of particles is a finite set $C$ of points together with a mass function $m:C\to\RR^+$. Since every point in an orbit of a symmetric configuration has the same mass, the mass function descends to a function $\bar m:(C/G)\to\RR^+$.}  

The following decomposition of $\CS(\Gamma)$ is immediate from the discussion above.

\begin{proposition}\label{prop:product}
Let $G<\OO(d)$ be a finite subgroup and $\Gamma=\sum a_H(H)$ a given Burnside type, where the sum is over $(H)\in\ISG$, and $a_H\in\NN$.  Then there is a diffeomorphism
$$\CS(\Gamma) \simeq \prod_{(H)\in\ISG}\; \CS(X_H,a_H),$$
where $\prod$ denotes the Cartesian product and $\CS(X_H,a)$ 
is the configuration space of $a$ particles in $X_H$.
\end{proposition}

Note that the connected components $X_H^\alpha$ of $X_H$ correspond to the topological Burnside types with orbit type $(H)$.  The expression above is readily refined to give 
\begin{equation} \label{eq:CS-compnts}
\CS(\Gamma) \simeq \prod_{(H)^\alpha\in\cT}\; \CS(X_H^\alpha,a_H^\alpha).
\end{equation}

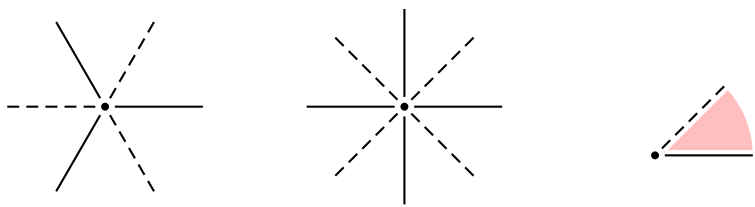
\begin{figure}[t] 
 \psset{unit=0.65}
 \centering
 \begin{pspicture}(-2,-2.5)(3,2)
   \psdot[dotsize=2pt](0,0)
    \psline(2,0)(0.2,0)
    \psline(-1,-1.732)(-0.1,-0.173) 
    \psline(-1,1.732)(-0.1,0.1732)
  \psset{linestyle=dashed}
    \psline(1,1.732)(.1,.1732)
    \psline(-2,0)(-0.2,0)
    \psline(1,-1.732)(.1,-.1732)
 \end{pspicture}
 \begin{pspicture}(-3,-2.5)(2,2)
   \psdot[dotsize=2pt](0,0)
    \psline(-2,0)(-0.2,0)
    \psline(0.2,0)(2,0)
    \psline(0,-2)(0,-0.2) 
    \psline(0,2)(0,0.2) 
  \psset{linestyle=dashed}
    \psline(1.414,1.414)(.1414,.1414)
    \psline(-1.414,1.414)(-.1414,.1414)
    \psline(1.414,-1.414)(.1414,-.1414)
    \psline(-1.414,-1.414)(-.1414,-.1414)
 \end{pspicture}
 \begin{pspicture}(-3,-1.5)(2,3)
   \pswedge*[linecolor=pink](0,0){2}{0}{45}
   \psdot[dotsize=5pt,linecolor=white](0,0)
   \psdot[dotsize=2pt](0,0)
    \psline[linewidth=4pt,linecolor=white](0.2,0)(2,0)
    \psline(0.2,0)(2,0)
    \psline[linewidth=4pt,linecolor=white](1.414,1.414)(.1414,.1414)
  \psset{linestyle=dashed}
    \psline(1.414,1.414)(.1414,.1414)
 \end{pspicture}
 \begin{minipage}{0.75\textwidth}
   \caption{The connected components of the fixed point sets for actions of $\D_3$ and $\D_4$, with the orbit space $\RR^2/\D_n$ with its 4 strata represented on the right.}
 \label{fig:dihedral quotient}
 \end{minipage}
 \end{figure}

One can identify each $X_H$ with a quotient of a subspace of $\RR^d$ as follows. For each isotropy subgroup $H$ let $V=V(H)=\Fix(H,\RR^d)$ and let $V^\circ=V^\circ(H)$ be the subset of points whose isotropy is precisely $H$ ($V^\circ$ is an open and dense subset of $V$). It follows from the relation $G_{g x} = gG_xg^{-1}$,  that for $x\in V^\circ$ one has $gx$ is also in $V^\circ$ if and only if $g\in N_G(H)$, the normaliser of $H$ in $G$.  It follows from this that $X_H\simeq V^\circ(H)/N_G(H)$. 
This is illustrated in Figure~\ref{fig:dihedral quotient} for the dihedral group $D_n$ acting on $\RR^2$.

We finish this section with some observations and an example.

\begin{itemize}
\item  
The trivial case where $V(H)=\{0\}$ does not contribute to the topology of $\CS(\Gamma)$. 

\item
The simplest non-trivial case is when the fixed point space is 1-dimensional: $\dim(V(H))=1$. This has already been mentioned in Section~\ref{sec:symmetry}, and is well-known. {If $a>1$ and the masses of the orbits of type $(H)$ are distinct then the space $C(X_H,a)$ is a disjoint union of contractible connected components, corresponding to different orderings of the points.  At the other extreme, if the masses are all equal then there is only a single component, since different orderings cannot be distinguished. See also Remark\,\ref{rmk:ordered}.}

\item 
If a fixed point space has dimension 2, there are two different possibilities. Let $H$ be the isotropy subgroup in question and $V$ the fixed point space (of dimension 2). The first possibility is that $V\setminus  V^\circ$ is a (finite) union of 1-dimensional subspaces, and in this case each component of $V^\circ$ is diffeomorphic to the plane. The contribution to $\CS(\Gamma)$, if $\Gamma$ includes $a(H)^\alpha$ is then diffeomorphic to $\CS(\RR^2,a)$.  {The topology of this space is well-known: its fundamental group is a subgroup of the braid group on $a$ strings depending on how many masses are equal (the pure braid group if they are all distinct or the full braid group if they are all equal),} while all its higher homotopy groups vanish \cite{FM-book}. 

The second possibility is that $V^\circ$ is a punctured plane, and the contribution to $\CS(\Gamma)$ from $a$ orbits of type $(H)$ is equivalent to $\CS(\RR^2,a+1)$. 

For example, consider the symmetric 5-body configurations, with symmetry $\D_1\simeq\ZZ_2$ acting by reflection in a line (see Figure\,2(c,d) of \cite{L-S09}). The orbit types are $(\D_1)$ and $(\one)$, and for a total of 5 bodies there are 3 possibilities,
$$1(\D_1)+2(\one),\quad 3(\D_1)+1(\one), \quad \mbox{and}\quad 5(\D_1).$$
The $5(\D_1)$ are the collinear Moulton configurations of 5 bodies.  For $\Gamma=1(\D_1) + 2(\one)$, the resulting space $\CS(\Gamma)$ is homotopic to the circle, so $U$ must have at least two critical points, as illustrated in \cite{L-S09}.  {As above, the topology in this last case in principle depends on whether the masses are equal or distinct; however, in this case of two particles in the plane, both spaces are homotopic to the circle. 
}

\item
Higher dimensional fixed points spaces will contribute to higher homotopy groups and cohomology, but the correspondence is not so easily understood. 

\item
Even though the set $(\RR^d)^\circ$ of points in $\RR^d$ with trivial isotropy may not be connected, the quotient $(\RR^d)^\circ/G$ is always connected. This is because the complement of $(\RR^d)^\circ$ is a union of linear subspaces, and so the only way $(\RR^d)^\circ$ is disconnected is through hyperplanes, and these only arise as fixed point sets for reflections (a matrix with $(d-1)$ eigenvalues equal to $+1$ and one equal to $(-1)$), and the reflection then identifies the two sides of the corresponding hyperplane. The following example shows that $(\RR^d)^\circ$  may not be contractible.
\end{itemize}

\begin{example}
  Consider finally Example \ref{eg:D2} above, and $\CS(\Gamma)$ for
  $\Gamma=\eps(G) + a(R_z) + b(R_y) + c(R_x) + e(\one)$. The orbit
  types $(G), (R_x),(R_y)$ and $(R_z)$ give spaces of dimension 1 or
  less, so if $a,b$ or $c>0$ their contribution is to increase the number of 
  connected components of $\CS(\Gamma)$, but not otherwise to change 
  its topology.  However the generic orbit, with orbit type $(\one)$, consists 
  of 4 points in the complement of the coordinate axes.  Its contribution
  to $\CS(\Gamma)$ is $\CS(X_\one,e)$, and one can show that the 
  stratum $X_\one=(\RR^3)^\circ/G$ is (homeomorphic to) the thrice 
  punctured sphere.  Thus each connected component of the space 
  $\CS(\Gamma)$ is homotopic to the space of $e$ points in the thrice 
  punctured sphere.  Using Morse theory one can show that for $e=1$, 
  and assuming critical points are non-degenerate, there must be at least 
  three critical points: one minimum and 2 saddle points, in each connected 
  component. In fact if $e=1$ and $a=b=c=0$ then there are 5 critical points: two 
  minima occurring at tetrahedral configurations and 3 saddles occurring 
  at squares in the coordinate hyperplanes.
\end{example}

\section{Balanced configurations} \label{sec:balanced}
Balanced configurations were introduced by Albouy and Chenciner \cite{A-C98} {as a configuration for which in a suitably larger space, the configuration is a relative equilibrium. More details are given in  \cite{Moeckel14}, and in \cite{Ch12} where several equivalent definitions are given. The version appropriate for our discussion is as follows.}  

Consider an ordered configuration $(x_1,\dots,x_n)$ in $(\RR^d)^n$ with $\sum m_ix_i=0$ and let $X$ be the $d\times n$ matrix whose columns are the position vectors of the points $x_1,\dots,x_n$, and let $\mu$ be the $n\times n$ diagonal matrix with $\mu_{ii}=m_i$. Consider the $d\times d$ matrix $S = X^T\mu X$.  It is clear that while $X$ depends on the order of the points, $S$ does not, so $S$ only depends on the configuration as defined in Section~\ref{sec:symmetry}. 

The \emph{inertia spectrum} of the configuration is the spectrum (with multiplicities) of the matrix $S$ \cite{Ch12}. It is easy to see that the eigenvalues are all non-negative, and that the moment of inertia function $I$ is equal to the trace $\tr(S)$. Moreover, 0 is an eigenvalue if and only if all the bodies are contained in a lower dimensional subspace. 

Let $\CS(\sigma)$ denote the space of all configurations with inertia spectrum $\sigma$.  A configuration is said to be \defn{balanced} if it is a critical point of the restriction of the potential function $U$ to $\CS(\sigma)$.  Since $I$ is constant on $\CS(\sigma)$ it follows that any central configuration is also a balanced configuration.  The variational argument used in this paper shows that there is always a balanced configuration in each non-empty $\CS(\sigma)$.

\begin{theorem}
Given any finite subgroup $G$ of\/ $\OO(d)$, let $\sigma$ be the inertia spectrum of some symmetric configuration .  Then $\CS(\sigma)^G$ is a non-empty closed subset of\/ $\CS(\sigma)$ and there is a symmetric balanced configuration in each component of\/ $\CS(\sigma)^G$, and indeed on $\CS(\sigma)(\Gamma)$, for any topological Burnside type $\Gamma$ for which $\CS(\sigma)(\Gamma)$ is non-empty. 
\end{theorem}

Symmetry of a configuration will cause its inertia spectrum to have multiplicities.  Given a finite subgroup $G<\OO(d)$,  decompose $\RR^d$ as a sum of isotypic representations of $G$, $\RR^d=\oplus_j E_J$.
That is, each $E_j$ is a sum of copies of isomorphic irreducible representations, and one can write $E_j=W_j\otimes\RR^{d_j}$, where $W_j$ is an irreducible representation and $d_j$ the multiplicity of that representation in $\RR^d$.  See for example the book of Serre \cite{Serre}.

It is clear that for a symmetric configuration, for each matrix $A\in G$, the symmetric matrix $S$ satisfies $A^TSA=S$. Since $A$ is orthogonal, $A^T=A^{-1}$ whence $AS=SA$ for all $A\in G$.  It then follows from Schur's Lemma \cite{Serre} that the matrix $S$ block diagonalizes into a single block $S_j$ for each $E_j$. On the $E_j$ block, the eigenvalues will have multiplicity at least $\dim W_j$. In the particular case that $d_j=1$, so $E_j=W_j$ is irreducible, the symmetric matrix $S_j$ will be a scalar matrix, equal to the moment of inertia of the projection of the configuration into $E_j$ times the identity.  Two immediate conclusions are as follows.

\begin{proposition}
Suppose $\RR^d$ is an irreducible representation of\/ $G<\OO(d)$. Then a symmetric configuration is balanced if and only if it is central.
\end{proposition}

For example, any balanced configuration with tetrahedral, octahedral or icosahedral symmetry is a central configuration (see Example~\ref{eg:Platonic}). 

{
\begin{proof}
In this case $S$ is a scalar matrix, and fixing the inertia spectrum $\sigma$ is equivalent to fixing the  moment of inertia $I$.
\end{proof}
}

\begin{proposition}
Suppose $\RR^d=\oplus_j E_j$ is a sum of distinct irreducible representations of $G<\OO(d)$,  (ie no representation occurs with multiplicity $>1$). Then two configurations $C_1,C_2$ with the same Burnside type have the same spectrum if and only if, for each $j$, their images $\pi_j(C_i)$ have the same moment of inertia, where $\pi_j:\RR^d\to E_j$ is the natural projection. 
\end{proposition}


\paragraph{Acknowledgements} I would like to thank Manuele Santoprete for pointing out some references, and Alain Albouy for making a number of helpful suggestions.

\let\oldbibliography\thebibliography
  \renewcommand{\thebibliography}[1]{%
  \oldbibliography{#1}
  \setlength{\itemsep}{-3pt}
  \small
  } 

\end{document}